\documentclass[12pt]{amsart}
\usepackage{amscd,amsmath,amsthm,amssymb,graphics}
\usepackage{pstcol,pst-plot,pst-3d}
\usepackage{lmodern,pst-node}
\usepackage{multicol}
\usepackage{amssymb}
\usepackage{graphicx}
\usepackage{tikz}
\usepackage{subfig}
\usepackage{mathtools}
\usepackage{tikz-3dplot}
\usepackage{enumitem}
\usepackage{graphicx}
\usepackage{tikz}
\usepackage{enumitem}
\usepackage{enumerate}
\usepackage{epic,eepic}
\usepackage{amsfonts,amssymb,amscd,amsmath,enumerate,verbatim}
\psset{unit=0.7cm,linewidth=0.8pt,arrowsize=2.5pt 4}

\newpsstyle{fatline}{linewidth=1.5pt}
\newpsstyle{fyp}{fillstyle=solid,fillcolor=verylight}
\definecolor{verylight}{gray}{0.97}
\definecolor{light}{gray}{0.9}
\definecolor{medium}{gray}{0.85}
\definecolor{dark}{gray}{0.6}

\def\NZQ{\Bbb}               

\def\RR{{\NZQ R}}

\def\frk{\frak}               

\def\Phi{{\frk n}}
\def\Phi{{\frk N}}
\def\MI{{\mathcal I}}

\def\MO{{\mathcal O}}
\def\MC{{\mathcal C}}
\def\xb{{\bold x}}

\def\eb{{\bold e}}
\def\opn#1#2{\def#1{\operatorname{#2}}} 
\opn\chara{char} \opn\length{\ell} \opn\pd{pd} \opn\rk{rk}
\opn\projdim{proj\,dim} \opn\injdim{inj\,dim} \opn\rank{rank}
\opn\depth{depth} \opn\grade{grade} \opn\height{height}
\opn\embdim{emb\,dim} \opn\codim{codim}

\opn\Tr{Tr} \opn\bigrank{big\,rank}
\opn\superheight{superheight}\opn\lcm{lcm}
\opn\trdeg{tr\,deg}
\opn\reg{reg} \opn\lreg{lreg} \opn\ini{in} \opn\lpd{lpd}
\opn\size{size}\opn\bigsize{bigsize}
\opn\cosize{cosize}\opn\bigcosize{bigcosize}
\opn\sdepth{sdepth}\opn\sreg{sreg}
\opn\link{link}\opn\fdepth{fdepth}
\opn\div{div} \opn\Div{Div} \opn\cl{cl} \opn\Cl{Cl}
\opn\Spec{Spec} \opn\Supp{Supp} \opn\supp{supp} \opn\Sing{Sing}
\opn\Ass{Ass} \opn\Min{Min}\opn\Mon{Mon} \opn\dstab{dstab} \opn\astab{astab}
\opn\Syz{Syz}
\opn\Ann{Ann} \opn\Rad{Rad} \opn\Soc{Soc}
\opn\Im{Im} \opn\Ker{Ker} \opn\Coker{Coker} \opn\Am{Am}
\opn\Hom{Hom} \opn\Tor{Tor} \opn\Ext{Ext} \opn\End{End}
\opn\Aut{Aut} \opn\id{id}

\opn\nat{nat}
\opn\pff{pf}
\opn\Pf{Pf} \opn\GL{GL} \opn\SL{SL} \opn\mod{mod} \opn\ord{ord}
\opn\Gin{Gin} \opn\Hilb{Hilb}\opn\sort{sort}
\opn\S{S}
\opn\aff{aff} \opn\con{conv} \opn\relint{relint} \opn\st{st}
\opn\lk{lk} \opn\cn{cn} \opn\core{core} \opn\vol{vol}
\opn\link{link} \opn\star{star}\opn\lex{lex}\opn\ini{in}
\opn\gr{gr}

\def\pot#1#2{#1[\kern-0.28ex[#2]\kern-0.28ex]}

\opn\dirlim{\underrightarrow{\lim}}
\opn\inivlim{\underleftarrow{\lim}}
\let\to=\rightarrow

\def\Implies{\ifmmode\Longrightarrow \else
        \unskip${}\Longrightarrow{}$\ignorespaces\fi}
\def\implies{\ifmmode\Rightarrow \else
        \unskip${}\Rightarrow{}$\ignorespaces\fi}
\def\iff{\ifmmode\Longleftrightarrow \else
        \unskip${}\Longleftrightarrow{}$\ignorespaces\fi}

\let\:=\colon
\newtheorem{Theorem}{Theorem}

 \newtheorem{Definition}[Theorem]{Definition}

\let\epsilon\varepsilon
\let\kappa=\varkappa
\def\qed{\ifhmode\textqed\fi
      \ifmmode\ifinner\quad\qedsymbol\else\dispqed\fi\fi}
\def\textqed{\unskip\nobreak\penalty50
       \hskip2em\hbox{}\nobreak\hfil\qedsymbol
       \parfillskip=0pt \finalhyphendemerits=0}
\def\dispqed{\rlap{\qquad\qedsymbol}}

\opn\dis{dis}
\def\pnt{{\raise0.5mm\hbox{\large\bf.}}}

\opn\Lex{Lex}

\begin{document}
 \title {Compatible algebras with straightening laws on distributive lattices}

 \author {Daniel B\u anaru and Viviana Ene}

\address{Daniel B\u anaru, Faculty of Mathematics and Computer Science,  University of Bucharest, Academiei 14,
 010014 Bucharest, Romania}\email{danielbanaru@yahoo.ro}

\address{Viviana Ene, Faculty of Mathematics and Computer Science, Ovidius University, Bd.\ Mamaia 124,
 900527 Constanta, Romania}\email{vivian@univ-ovidius.ro}


 \begin{abstract} We characterize the finite distributive lattices on which there exists a unique compatible algebra 
with straightening laws.
 \end{abstract}

\thanks{}

 \maketitle

\section*{Introduction}

Let $P$ be  a finite partial ordered set (poset for short) and $\MI(P)$ the distributive lattice of the poset ideals of $P.$ By a famous 
theorem of Birkhoff \cite{B}, for every finite distributive lattice $L$, there exists a unique subposet $P$ of $L$ such that $L=\MI(P).$ 
The order polytope $\MO(P)$ and the chain polytope $\MC(P)$ were introduced in \cite{Sta}. In \cite{Hibi} it was shown that the toric ring 
$K[\MO(P)]$ over a field $K$ is an algebra with straightening laws (ASL in brief) on the distributive lattice $\MI(P)$ over the field $K.$ In \cite{HibiLi} it was shown that the ring $K[\MC(P)] $ associated with the chain polytope shares the same property. 

Let $S=K[x_1,\ldots,x_n,t]$ be the polynomial ring over a field $K$ and $\{w_\alpha\}_{\alpha\in \MI(P)}$ be a set of monomials in $x_1,\ldots,x_n$ indexed by $\MI(P)$. Let $K[\Omega]\subset S$ be the toric ring generated over $K$ by the set of  monomials 
$\Omega=\{\omega_\alpha\}_{\alpha\in \MI(P)}$ where $\omega_\alpha=w_\alpha t$ for all $\alpha\in \MI(P).$ Clearly, $K[\Omega]$ is a graded algebra if we set $\deg(\omega_\alpha)=1$ for all $\alpha \in \MI(P).$ Let $\varphi: \MI(P)\to K[\Omega]$ be the injective map defined by 
$\varphi(\alpha) =\omega_\alpha $ for all $\alpha\in \MI(P).$ Assume that $K[\Omega]$ is an ASL on $\MI(P)$ over $K.$ According to 
\cite{HibiLi}, $K[\Omega]$ is  a \textit{compatible } ASL if each of its straightening relations  is of the form 
$\varphi(\alpha)\varphi(\alpha^\prime)=\varphi(\beta)\varphi(\beta^{\prime})$ with $\beta\subseteq \alpha\cap\alpha^\prime$
and $\beta^\prime\supseteq \alpha\cup\alpha^\prime,$ where $\alpha,\alpha^\prime$ are incomparable elements in $\MI(P).$ If $K[\Omega]$ and 
$K[\Omega^\prime]$ are compatible ASL on $\MI(P)$ over $K,$ we identify them if they have the same straightening relations. In this case, we write $K[\Omega]\equiv K[\Omega^\prime].$

In \cite[Question 5.1 (b)]{HibiLi}, Hibi and Li asked the following question:\textit{ For which posets $P$, does there exists a unique 
compatible ASL on $\MI(P)$ over $K?$}

In this note, we give a complete answer to this question. Namely, we prove the following.

\begin{Theorem}\label{main}
Let $P$ be a finite poset. Then, the following statements are equivalent:
\begin{itemize}
	\item [{\em (i)}] There exists a unique  compatible ASL on $\MI(P)$ over $K.$
	\item[{\em (ii)}] $K[\MO(P)]\equiv K[\MC(P)]\equiv K[\MC(P^\ast)],$ where $P^\ast$ denotes the dual poset of $P.$
	\item[{\em (iii)}] Each connected component of $P$ is a chain, that is, $P$ is a direct sum of chains.
\end{itemize}
\end{Theorem}

\section{Order polytopes, chain polytopes, and their associated toric rings}

Let $P=\{p_1,\ldots,p_n\}$ be a finite poset. For the basic terminology regarding posets which is used in this paper we refer to 
\cite{B} and \cite[Chapter 3]{Sta2}. The order polytope $\MO(P)$ is defined as 
\[
\MO(P)=\{\xb=(x_1,\ldots,x_n)\in \RR^n: 0\leq x_i\leq 1,  1\leq i\leq n,  x_i\geq x_j \text{ if } p_i\leq p_j \text{ in }P\}.
\]

In \cite[Corollay 1.3]{Sta} it was shown that the vertices of $\MO(P)$ are $\sum_{p_i\in \alpha}\eb_i, \alpha\in\MI(P).$ Here 
$\eb_i$ denotes the unit coordinate vector in $\RR^n.$ If $\alpha=\emptyset,$ then the corresponding vertex in $\MO(P)$ is the 
origin of $\RR^n.$

The chain polytope $\MC(P)$ is defined as
\[
\MC(P)=\{\xb=(x_1,\ldots,x_n)\in \RR^n: x_i\geq 0, 1\leq i\leq n, \]\[\text{ and }x_{i_1}+\cdots+x_{i_r}\leq 1 \text{ if } p_{i_1}<\cdots <p_{i_r}
\text{is a maximal chain in } P\}.
\]
In \cite[Theorem 2.2]{Sta}, it was proved that the vertices of $\MC(P)$ are $\sum_{p_i\in A}\eb_i$ where $A$ is an antichain in $P.$ Recall that an antichain in $P$ is a subset of $P$ such that any two distinct elements in the subset are incomparable. Since every poset ideal is uniquely determined by its antichain of maximal elements, it follows that $\MO(P)$ and $\MC(P)$ have the same number of vertices. However, as it was observed in \cite{Sta}, in general, these two polytopes are not combinatorially equivalent, that is, 
$\MO(P)$ and $\MC(P)$ need not have the same number of $i$-dimensional faces for $i>0.$ Combinatorially equivalence of order and chain polytopes is studied in \cite{HibiLi2}.

\subsection*{The toric rings $K[\MO(P)]$ and $K[\MC(P)]$} To each subset $W\subset P$ we attach the squarefree monomial 
$u_W\in K[x_1,\ldots,x_n],$ $u_W=\prod_{p_i\in W}x_i.$ If $W=\emptyset,$ then $u_W=1.$ The toric ring $K[\MO(P)]$, known as the Hibi ring associated with  the distributive lattice $\MI(P),$ is generated over $K$ by all the monomials $u_\alpha t\in S$ where $\alpha\in \MI(P).$
The toric ring $K[\MC(P)]$ is generated by all the monomials $u_A t$ where $A$ is an antichain in $P.$  Also, as we have already mentioned in Introduction, both rings are algebras with straightening laws on $\MI(P)$ over $K.$

We recall the definition of an ASL. For a quick introduction to this topic we refer to \cite{Eis} and \cite[Chapter XIII]{Hibi2}.  Let $K$ be a field,  $R=\oplus_{i\geq 0}R_i$ with $R_0=K$ be a graded $K$-algebra, $H$ a finite poset, and $\varphi:H\to R$ an injective map which maps each $\alpha\in H$ to a homogeneous element $\varphi(\alpha)\in R$ with 
$\deg \varphi(\alpha)\geq 1$. A \textit{standard monomial} in $R$ is a product  
$\varphi(\alpha_1)\varphi(\alpha_2)\cdots \varphi(\alpha_k)$ where $\alpha_1\leq \alpha_2\leq \cdots \leq \alpha_k$ in $H.$

\begin{Definition}\label{defasl}{\em 
The $K$-algebra $R$ is called an algebra with straightening laws on $H$ over $K$ if the following conditions hold:
\begin{itemize}
	\item[(1)] The set of standard monomials is a $K$--basis of $R;$
	\item[(2)] If $\alpha,\beta\in H$ are incomparable and if $\alpha\beta=\sum c_i \gamma_{i1}\ldots \gamma_{ik_i},$ where $c_i\in K\setminus\{0\}$ and $\gamma_{i1}\leq \ldots \leq \gamma_{ik_i},$ is the unique expression of $\alpha\beta$ as a linear combination of standard monomials, then $\gamma_{i1}\leq \alpha, \beta$ for all $i.$
\end{itemize}}
\end{Definition}
The above  relations $\alpha\beta=\sum c_i \gamma_{i1}\ldots \gamma_{ik_i}$ are called the straightening relations of $R$ and they generate the defining ideal of $R.$

Let us go back to the toric rings $K[\MO(P)]$ and $K[\MC(P)].$ 

One considers $\varphi: \MI(P)\to K[\MO(P)]$ defined by 
$\varphi(\alpha)=u_\alpha t$ for every $\alpha\in \MI(P).$ As it was proved by Hibi in \cite{Hibi}, $K[\MO(P)]$ is an 
ASL on $\MI(P)$ over $K$ with the straightening relations $\varphi(\alpha)\varphi(\beta)=\varphi(\alpha\cap\beta)\varphi(\alpha\cup\beta)$
where $\alpha,\beta$ are incomparable elements in $\MI(P).$

On the other hand, one defines $\psi:\MI(P)\to K[\MC(P)]$ by setting $\psi(\alpha)=u_{\max\alpha}t$ for all $\alpha\in \MI(P)$ where 
$\max \alpha$ denotes the set of the maximal elements in $\alpha.$ Note that, for every $\alpha\in \MI(P),$ $\max \alpha$ is an antichain 
in $P$ and each antichain $A\subset P$ determines a unique ideal $\alpha\in \MI(P),$ namely, the poset ideal generated by $A.$ Therefore, 
$\psi$ is an injective well defined map and by \cite[Theorem 3.1]{HibiLi}, the ring $K[\MC(P)]$ is an ASL on $\MI(P)$ over $K$ with the 
straightening relations \[\psi(\alpha)\psi(\beta)=\psi(\alpha\ast\beta)\psi(\alpha\cup\beta)\] where $\alpha\ast\beta$ is the poset ideal of $P$ generated by $\max(\alpha\cap\beta)\cap(\max\alpha\cup\max\beta).$

We observe that one may also consider $K[\MC(P^\ast)]$ as an  ASL on $\MI(P)$, where $P^\ast$ is the dual poset of $P.$ We may define
$\delta:\MI(P)\to K[C(P^\ast)]$ by $\delta(\alpha)=u_{\min\overline{\alpha}}t$ for $\alpha\in \MI(P),$ where $\min\overline{\alpha}$ is the set of minimal elements in $\overline{\alpha}$ and $\overline{\alpha}$ is the filter  $P\setminus \alpha$ of $P.$ We recall that a 
\textit{filter} $\gamma$ in $P$ (or \textit{dual order ideal}) is a subset of $P$ with the property that for every $p\in \gamma$ and every $q\in P$ with $q\geq p,$ we have $q\in \gamma.$ Thus, a filter in $P$ is simply a poset ideal in the dual poset $P^\ast.$ The ring $K[C(P^\ast)]$ is an ASL on $\MI(P)$ over $K$ as well with the straightening relations \[\delta(\alpha)\delta(\beta)=\delta(\alpha\cap\beta)\delta(\alpha\circ\beta)\] for 
incomparable elements $\alpha,\beta\in \MI(P),$ where $\alpha\circ \beta$ is the poset ideal of $P$ which is the complement in $P$ of the filter generated by 
$\min(\overline{\alpha}\cap\overline{\beta})\cap(\min\overline{\alpha} \cup \min\overline{\beta}).$ Let us also observe that all the algebras $K[\MO(P)], K[\MC(P)]$, and 
$K[\MC(P^\ast)]$ are compatible algebras with straightening laws.

\section{Proof of Theorem~\ref{main}}

We clearly have (i) $\Rightarrow$ (ii). Let us now prove (ii) $\Rightarrow$ (iii). By hypothesis, the straightening relations of 
$K[\MO(P)], K[\MC(P)],$ and $K[\MC(P^\ast)]$ coincide. Therefore, we must have 
\begin{equation}\label{eq1}
\alpha\cap \beta=\alpha \ast \beta \text{ and }
\overline{\alpha}\cap\overline{\beta}=\overline{\alpha\circ \beta}
\end{equation}
for all $\alpha,\beta$ incomparable elements in $\MI(P).$ From the second equality in (\ref{eq1}), it follows that  $\overline{\alpha}\cap\overline{\beta}$ is the filter of $P$ generated by $\min(\overline{\alpha}\cap\overline{\beta})\cap(\min\overline{\alpha} \cup \min\overline{\beta}).$
Assume that there exists two incomparable elements $p,p^\prime\in P$ such that there exists $ q\in P$ with $q>p$ and $q>p^\prime.$ Consider $\overline{\alpha}$ the filter generated by $p$ and $\overline{\beta}$ the filter generated by $p^\prime.$ Then 
$\min(\overline{\alpha}\cap\overline{\beta})\cap(\min\overline{\alpha} \cup \min\overline{\beta})=\emptyset$, but obviously, 
$ \overline{\alpha}\cap\overline{\beta}\neq \emptyset.$ This shows that, for any two incomparable elements $p,p^\prime\in P,$  there is no element $q\in P$ with $q>p, q>p^\prime.$

Similarly, by using the first equality in (\ref{eq1}), we derive that, for any two incomparable elements $p,p^\prime\in P,$ there is no element $q\in P$ with $q<p, q<p^\prime$. This shows that every connected component of the poset $P$ is a chain.

Finally, we prove (iii) $\Rightarrow$ (i). Let $P$ be a poset  such that all its connected components are chains and assume that the cardinality of $P$ is equal to $n.$ Let $\{\omega_\alpha\}_{\alpha\in \MI(P)}$ be the generators of $K[\Omega]\subset S$ and assume that the straightening relations 
of $K[\Omega]$ are $\varphi(\alpha)\varphi(\alpha^\prime)=\varphi(\beta)\varphi(\beta^\prime)$ where 
$\beta\subseteq \alpha\cap\alpha^\prime,$  $\beta^\prime\supseteq \alpha\cup\alpha^\prime,$ and $\alpha,\alpha^\prime$ are incomparable elements in $\MI(P).$ We have to show that, for all $\alpha,\alpha^\prime$  incomparable elements in $\MI(P),$ we have 
$\beta= \alpha\cap\alpha^\prime$ and $\beta^\prime= \alpha\cup\alpha^\prime.$ 

We proceed by induction on 
\[
k=n-(\rank(\alpha\cup \alpha^\prime)-\rank(\alpha\cap \alpha^\prime)).
\]
If $k=0,$ that is, $\rank(\alpha\cup \alpha^\prime)-\rank(\alpha\cap \alpha^\prime)=n$, then $\alpha\cup \alpha^\prime=P$ and 
$\alpha\cap \alpha^\prime=\emptyset,$ thus $\beta= \alpha\cap\alpha^\prime$ and $\beta^\prime= \alpha\cup\alpha^\prime.$ Assume that the desired conclusion is true for $\rank(\alpha\cup \alpha^\prime)-\rank(\alpha\cap \alpha^\prime)=n-k$ with $k\geq 0.$ Let us choose now 
$\alpha,\alpha^\prime$ incomparable in $\MI(P)$ such that $\rank(\alpha\cup \alpha^\prime)-\rank(\alpha\cap \alpha^\prime)=n-k-1$ 
and assume that we have a straightening relation $\varphi(\alpha)\varphi(\alpha^\prime)=\varphi(\beta)\varphi(\beta^\prime)$ with 
$\beta\subsetneq \alpha\cap\alpha^\prime$ or $\beta^\prime\supsetneq \alpha\cup\alpha^\prime.$ By duality, we may reduce to considering  
$\beta^\prime\supsetneq \alpha\cup\alpha^\prime.$ In other words, in $K[\Omega]$ we have 
\[
\omega_{\alpha}\omega_{\alpha^\prime}=\omega_{\beta}\omega_{\beta^\prime}, \text{ with } \beta\subseteq \alpha\cap\alpha^\prime
\text{ and }\beta^\prime\supsetneq \alpha\cup\alpha^\prime.
\] As $P$ is a direct sum of chains, we may find $p\in \max(\alpha\cup\alpha^\prime)$ and  $q\in \beta^\prime\setminus(\alpha\cup\alpha^\prime)$
 such that $q$ covers $p$ in $P$, 
that is, there is no other element $q^\prime$ in $P$ with $q>q^\prime >p.$ Without loss of generality, we may assume that $p\in \alpha^\prime.$ Let $\alpha_1$ be the poset ideal of $P$ generated by $\alpha^\prime\cup\{q\}.$ As all the connected components of $P$ are chains, we have $\alpha_1=\alpha^\prime\cup \{q\}$ since there are no other elements in $P$ which are smaller than $q$ except those which are on the same chain as $p$ and $q$ which are  in $\alpha^\prime.$ Moreover, by the choice of $q,$ we have 
\[\alpha_1\subseteq \beta^\prime  \text{ and }  \alpha\cap\alpha_1=\alpha\cap(\alpha^\prime\cup\{q\})=\alpha\cap\alpha^\prime.\]
On the other hand, 
\[
\rank(\alpha\cup\alpha_1)-\rank(\alpha\cap\alpha_1)=\rank(\alpha\cup\alpha^\prime\cup\{q\})-\rank(\alpha\cap\alpha^\prime)=\]
\[=\rank(\alpha\cup\alpha^\prime)+1-\rank(\alpha\cap\alpha^\prime)=n-k.
\]
By the inductive hypothesis, it follows that $\varphi(\alpha)\varphi(\alpha_1)=\varphi(\alpha\cap\alpha_1)\varphi(\alpha\cup\alpha_1),$ or, 
equivalently, in $K[\Omega]$ we have the equality 
$\omega_\alpha \omega_{\alpha_1}=\omega_{\alpha\cap\alpha_1}\omega_{\alpha\cup\alpha_1}.$ Thus, we have obtained the following equalities 
in $K[\Omega]:$
\[
\omega_{\alpha}\omega_{\alpha^\prime}=\omega_{\beta}\omega_{\beta^\prime} \text{ and }\omega_\alpha \omega_{\alpha_1}=
\omega_{\alpha\cap\alpha^\prime}\omega_{\alpha\cup\alpha_1}.\] This implies that 
\begin{equation}\label{eq3}
\omega_{\alpha\cap\alpha^\prime}\omega_{\alpha^\prime}\omega_{\alpha\cup \alpha_1}=\omega_\beta \omega_{\alpha_1}\omega_{\beta^\prime}.
\end{equation}
In addition, we have:
 \[\alpha\cap\alpha^\prime\subset \alpha^\prime\subset \alpha\cup\alpha_1 \text{ and }  \beta\subset \alpha\cap\alpha^\prime=\alpha\cap\alpha_1\subset \alpha_1\subset \beta^\prime.\]
This implies that both monomials in (\ref{eq3}) are standard monomials in $K[\Omega]$ which is in contradiction to the condition that the standard monomials form a $K$-nasis in $K[\Omega].$ Therefore, our proof is completed.

\end{document}